\documentclass[psamsfonts]{amsart}

\usepackage{amssymb,amsfonts}
\usepackage[all,arc]{xy}
\usepackage{enumerate}
\usepackage{mathrsfs}
\usepackage{graphicx}

\newtheorem{theorem}{Theorem}[section]
\newtheorem{corollary}[theorem]{Corollary}
\newtheorem{proposition}[theorem]{Proposition}
\newtheorem{lemma}[theorem]{Lemma}

\theoremstyle{definition}
\newtheorem{definition}[theorem]{Definition}

\newtheorem{example}[theorem]{Example}
\newtheorem{examples}[theorem]{Examples}

\newtheorem{remark}[theorem]{Remark}

\DeclareMathOperator{\Nil}{Nil}
\DeclareMathOperator{\Id}{Id}
\DeclareMathOperator{\Spec}{Spec}

\DeclareMathOperator{\Ann}{Ann}
\DeclareMathOperator{\Ass}{Ass}
\DeclareMathOperator{\ord}{ord}


\numberwithin{equation}{section}

\begin{document}

\title[Distinguished elements in semiring extensions]{Distinguished elements in semiring extensions}

\author[Peyman Nasehpour]{\bfseries Peyman Nasehpour}

\email{nasehpour@gmail.com}

\dedicatory{In memory of my father Maestro Nasrollah Nasehpour (1940 -- 2023)}



\begin{abstract}
In this paper, we investigate zero-divisor, nilpotent, idempotent, unit, small, and irreducible elements in semiring extensions such as amount, content, and monoid semialgebras. We also introduce new concepts such as the prime avoidance property in semirings, entire-like semirings, semialgebras with Property (A), and also, Armendariz and McCoy semialgebras and we prove some results related to these concepts. For example, we prove that if $B$ is an $S$-semialgebra, then under some conditions, the set of zero-divisors $Z(B)$ of $B$ is the union of the extended maximal primes of $Z(S)$. Finally, we prove a generalization of Eisenstein's irreducibility criterion.
\end{abstract}

\maketitle

\section{Introduction}

The main purpose of this paper is to investigate zero-divisor, nilpotent, idempotent, invertible, and other distinguished elements in semiring extensions. Though some of the results in the current paper are generalizations of their counterparts in ring theory, the other results are exclusively correct for a proper semiring due to the nature of semirings (because semiring is a kind of ring-like algebraic structure in which subtraction either does not exists or is disallowed).   

Since the language for semirings is not standardized \cite{Glazek2002}, we introduce a couple of concepts and fix some terminologies. In this paper, by a semiring \cite{Golan1999(a), Golan1999(b), Golan2003}, we mean a nonempty set $S$ with two binary operations of addition and multiplication such that the following conditions are satisfied:

\begin{enumerate}
	\item $(S,+)$ is a commutative monoid with identity element $0$;
	\item $(S,\cdot)$ is a commutative monoid with identity element $1 \not= 0$;
	\item Multiplication distributes over addition, i.e., $a\cdot (b+c) = a \cdot b + a \cdot c$ for all $a,b,c \in S$;
	\item The element $0$ is the absorbing element of the multiplication, i.e., $s \cdot 0=0$ for all $s\in S$.
\end{enumerate}

A semiring is proper if it is not a ring. Let $S$ be a semiring and $I$ a nonempty subset of $S$. The set $I$ is said to be an ideal of $S$, if $a+b \in I$ for all $a,b \in I$ and $sa \in I$ for all $s \in S$ and $a \in I$ (introduced for the first time by Bourne in \cite{Bourne1951}). We collect all ideals of a semiring $S$ in $\Id(S)$. The ideal $I$ of $S$ is called proper if $I \neq S$. A proper ideal $\mathfrak{p}$ of $S$ is prime, if $ab\in \mathfrak{p}$ implies either $a\in \mathfrak{p}$ or $b\in \mathfrak{p}$. We collect all prime ideals of a semiring $S$ in $\Spec(S)$. Finally, an ideal $I$ is subtractive if $a+b\in I$ and $a\in I$ imply that $b\in I$ for all $a,b\in S$ \cite{Golan1999(b)}.

By a semiring morphism, it is meant a function $\gamma$ from a semiring $S$ into a semiring $T$ satisfying the following conditions:

\begin{enumerate}
	\item $\gamma(a+b) = \gamma(a) + \gamma(b)$ and $\gamma(ab) = \gamma(a) \gamma(b)$, for all $a$ and $b$ in $S$.
	\item $\gamma(0) = 0$ and $\gamma(1) = 1$.
\end{enumerate}

In this paper, an extension of a semiring $S$ is a semiring $B$ such that there is a semiring morphism from $S$ into $B$. In this case, we say that $B$ is an $S$-semialgebra.

In the first section of the paper, we introduce amount functions over semialgebras that are generalizations of content functions over semiring polynomials and formal power series. Then, we investigate those basic properties of amount functions that we need in the other sections of the current paper. 

Let $S$ be a semiring and $B$ an $S$-semialgebra. We say a function $A$ from $B$ to the set of ideals $\Id(S)$ of $S$, defined by $f \mapsto A_f$, is an amount function if the following properties hold for all $s\in S$ and $f,g\in B$:

\begin{enumerate}
	\item $A$ preserves 0 and 1, i.e. $A_0 = (0)$ and $A_1 = S$.
	\item If $A_f =(0)$, then $f=0$. 
	\item $A$ is homogeneous, i.e. $A_{sf} = sA_f$.
	\item $A$ is submultiplicative, i.e. $A_{fg} \subseteq A_f A_g$.
\end{enumerate}

Note that in this paper, we use ``$\subseteq$'' for inclusion (see Definition 1.10 in \cite{Monk1969}).

Let us recall that if $f$ is a polynomial or power series over $S$, its content, denoted by $A_f$ \cite{GilmerGramsParker1975} or $c(f)$ \cite{Northcott1959}, is an ideal of $S$ generated by its coefficients. It is easy to check that such functions satisfy the conditions of amount functions mentioned above. 

In \S\ref{sec:amountsemialgebras}, we also introduce amount semialgebras (see Definition \ref{amountsemialgebras}) and in Theorem \ref{amountisweakcontent}, prove that if $B$ is an amount $S$-semialgebra, then $A_f A_g \subseteq \sqrt{A_{fg}}$ for all $f,g\in B$.

Let us recall that a commutative ring $R$ has Property (A) if for any finitely generated ideal $I$ of $R$ the inclusion $I \subseteq Z(R)$ implies that $I$ has a nonzero annihilator. As a generalization to this concept, we introduce the concept of a semialgebra having Property (A) (check Definition \ref{PropertyAsemialgebradef}) and prove that if the semiring $S$ has very few zero-divisors, then the $S$-semialgebra $S[X]$ has Property (A) (see Corollary \ref{PropertyAsemialgebracor}). Note that a semiring $S$, by definition, has very few zero-divisors if the set of zero-divisors $Z(S)$ of $S$ is a finite union of prime ideals in $\Ass(S)$ (see Definition 48 in \cite{Nasehpour2016}). The concept of rings and modules having very few zero-divisors has been investigated in \cite{NasehpourPayrovi2010,Nasehpour2010,Nasehpour2011}. 

Let $S$ be a semiring and $\Gamma \subseteq \Spec(S)$. We say $Z(S)$ is controlled by $\Gamma$ if $Z(S) = \bigcup_{\mathfrak{p} \in \Gamma} \mathfrak{p}$ (see Definition \ref{controllingzerodivisors}). In \S\ref{sec:PropertyAsemialgebras}, we also introduce the ``prime avoidance'' property for ideals of semirings which is a counterpart of its concept in commutative ring theory introduced by Chen in \cite{Chen2021}. Then, in Theorem \ref{PropertyAsemialgebrathm}, we prove that if $S$ is a semiring, $\Gamma \subseteq \Ass(S)$ satisfies the prime avoidance property, and $Z(S)$ is controlled by $\Gamma$, then any amount $S$-semialgebra has Property (A). This gives plenty of examples of semialgebras having Property (A) (check Examples \ref{PrimeAvoidanceExamples} and Examples \ref{ExamplesAmountSemialgebras}).

In \S\ref{sec:armendarizmccoysemialgebras}, we introduce Armendariz and McCoy semialgebras with respect to amount functions. Let $A$ be an amount function from the $S$-semialgebra $B$ into $\Id(S)$. We say $B$ is an Armendariz $S$-semialgebra if $fg = 0$ implies $A_f A_g = (0)$, for all $f$ and $g$ in $B$. In Proposition \ref{AmountisArmendariz}, we show that any amount semialgebra over a nilpotent-free semiring is an Armendariz semialgebra. Let us recall that a semiring $S$ is nilpotent-free if $s^2 = 0$ implies that $s = 0$, for all $s\in S$. As a generalization to Proposition 3.5 in \cite{GilmerGramsParker1975}, we show that if $S$ is a semiring such that each minimal prime ideal of $S$ is subtractive and  $f$ and $g$ are elements of $S[[X]]$ such that $fg = 0$, then $A_f A_g \subseteq \Nil(S)$ (see Proposition \ref{proposition35GGP}). A corollary to this is that if $S$ is nilpotent-free and each minimal prime ideal of $S$ is subtractive, then $S[[X]]$ is an Armendariz $S$-semialgebra. Note that semirings in which all prime ideals are subtractive (i.e. weak Gaussian semirings) have been investigated in \S3 in \cite{Nasehpour2016}.

McCoy \cite{McCoy1942} proved that if $R$ is a commutative ring and $fg = 0$ for some $f$ and $g$ in $R[X]$ with $g \neq 0$, then there is a nonzero $r\in R$ such that $fr = 0$. Note that the equality $fr = 0$ is equivalent to $r c(f) = (0)$. As a generalization to McCoy's Theorem, the author has showed that if $G$ is a cancellative torsion-free commutative monoid, then the $S$-semialgebra $S[G]$ satisfies the McCoy's property (see Corollary 2.2 in \cite{Nasehpour2021}). Based on this, we define that if $B$ is an $S$-semialgebra and $A: B \rightarrow \Id(S)$ an amount function, then $B$ is a McCoy $S$-semialgebra if for all $f,g \in B$, $fg = 0$ together with $g \neq 0$ implies that there is a nonzero $s \in S$ such that $s A_f = (0)$ (see Definition \ref{McCoysemialgebradef}). In this direction, in Theorem \ref{powerseriesismccoy}, we prove that the $S$-semialgebra $S[[X]]$ is a McCoy $S$-semialgebra if one of the following conditions hold:

\begin{enumerate}
	\item The semiring $S$ is nilpotent-free and weak Gaussian.
	
	\item The semiring $S$ is zerosumfree (note that a semiring $S$ is zerosumfree if $a+b = 0$ implies that $a=b=0$, for all $a$ and $b$ in $S$ \cite[p. 7]{Golan2003}).
\end{enumerate}

At the end of \S\ref{sec:armendarizmccoysemialgebras}, we find some other examples for Armendariz semialgebras. Namely in Theorem \ref{moreexamplesarmendarizzemialgebras}, we show that if $S$ is a zerosumfree semiring, any monoid semialgebra and formal power series over $S$ are Armendariz. Here is a good place to emphasize that in this paper all monoids are commutative unless otherwise explicitly stated. 

We devote \S\ref{sec:zerodivisorsnilpotents} to zero-divisors and nilpotent elements. For example, in Theorem \ref{extensionofzerodivisors1}, we prove that if $S$ is a semiring, $\Gamma \subseteq \Ass(S)$ satisfies the prime avoidance property, $B$ is an amount $S$-semialgebra, where $A$ is its amount function from $B$ into $\Id(S)$, and if $Z(S)$ is controlled by $\Gamma$ and $B$ is a McCoy $S$-semialgebra, then $Z(B)$ is controlled by $\mathfrak{p}^{\epsilon}$s, where $\mathfrak{p} \in \Gamma$. On the other hand, in Theorem \ref{extensionofzerodivisors2}, we show that if $S$ is a semiring and $\Gamma \subseteq \Spec(S)$ satisfies the prime avoidance property, also, if $B$ is an amount $S$-semialgebra, where $A$ is its amount function from $B$ into $\Id(S)$, $Z(S)$ is controlled by $\Gamma$ and the $S$-semialgebra $B$ is McCoy and has Property (A), then $Z(B)$ is controlled by $\mathfrak{p}^{\epsilon}$s, where $\mathfrak{p} \in \Gamma$.

In \S\ref{sec:zerodivisorsnilpotents}, we also investigate the nilpotent elements in some semiring extensions (recall that the set of all nilpotent elements of $S$ is denoted by $\Nil(S)$). For example, in Theorem \ref{nilpotentsmonoidsemirings}, we prove that if $M$ is a monoid and $S$ a semiring, then $\Nil(S[M]) = \Nil(S)[M]$ if one of the following conditions holds:

\begin{enumerate}
	\item $S$ is a zerosumfree semiring.
	\item $S$ is a weak Gaussian semiring and $M$ is a totally ordered monoid.
\end{enumerate}

It is evident that any nilpotent element of a semiring $S$ is a zero-divisor on $S$. In Definition \ref{entirelikesemiringsdef}, we define a semiring $S$ to be entire-like if $Z(S) \subseteq \Nil(S)$. Then, in Proposition \ref{monoidsemiringentirelikepro}, we show that if $S$ is an entire-like semiring and the monoid semiring $S[M]$ is a McCoy $S$-semialgebra, then $S[M]$ is entire-like. The rest of \S\ref{sec:zerodivisorsnilpotents} is devoted to entire semirings, i.e., semirings without nontrivial zero-divisors. For example, in Theorem \ref{monoidsemiringzerosumfree}, we prove that if $S$ is a semiring and $M$ is a totally ordered monoid, then the following statements are equivalent:
	
\begin{enumerate}
\item For all $f$ and $g$ in $S[M]$, we have $\deg(f+g) = \max \{\deg(f),\deg(g)\}$.
\item The semiring $S$ is zerosumfree.
\item The semiring $S[M]$ is zerosumfree.
\end{enumerate}

Note that, in this paper, $(M,+,0,<)$ is a totally ordered monoid if there is a total order $<$ on the monoid $M$ compatible with its binary operation, i.e. $a < b$ implies $a+c < b+c$ for all $a,b$ and $c$ in $M$. We say $(M,+,0,<)$ is a positive monoid if $(M,+,0,<)$ is a totally ordered monoid and its neutral element 0 is its least element.

In \S\ref{sec:idempotentunitssmall}, we investigate other distinguished elements of semiring extensions, namely, idempotents, units, and small elements. For instance, in Proposition \ref{exampleskpotent}, we show that if $\{m_i\}_{i=1}^n$ is a subgroup of the monoid $M$ and $n\cdot 1=\sum_{i=1}^{n} 1$ is a unit in a semiring $S$, then $\displaystyle \frac{\sum_{i=1}^{n} X^{m_i}}{n\cdot 1}$ is a multiplicatively idempotent element of $S[M]$. Examples of semirings in which the elements of the form $n \cdot 1$ are invertible ($n \in \mathbb{N}$) are Gel'fand semirings since a semiring $S$ is Gel'fand if $1+s \in U(S)$ for all $s\in S$. \cite[p. 56]{Golan1999(b)}. We recall that an element $s$ of a semiring $S$ is a unit if there is an element $t$ in $S$ such that $st=1$ \cite[p. 50]{Golan1999(b)}. We denote the set of all units of a semiring $S$ by $U(S)$. Note that $U(S)$ is always nonempty since $1\in U(S)$.

Let us recall that a semiring $S$ is additively cancellative if $a+b = a+c$ implies $b=c$ for all $a,b,c \in S$ \cite[p. 49]{Golan1999(b)}. In Theorem \ref{idempotentacsemiring}, we prove that if $S$ is an additively cancellative semiring and $X_i$ is an indeterminate over $S$ for each positive integer $i \leq k$, then multiplicatively idempotent elements of $S[[X_1,X_2,\ldots,X_k]]$ and $S$ are the same. Note that Theorem \ref{idempotentacsemiring} is a semiring version of Proposition 2 in \cite{KanwarLeroyMatczuk2013}. However, if a semiring $S$ is not additively cancellative, then Theorem \ref{idempotentacsemiring} may fail as we give some examples (see Examples \ref{idempotentacsemiringex}) to verify our claim.

As mentioned above, we also investigate unit and small elements of semialgebras in \S\ref{sec:idempotentunitssmall}. For example, we show that if $S$ is zerosumfree, then \[U(S[[X_1,X_2,\ldots,X_n]]) = U(S),\] where $X_1,X_2,\ldots,X_n$ are indeterminates over $S$ (check Theorem \ref{unitspowerseriesoverzerosumfreethm}).

Let us recall that an element $a$ of a semiring $S$ is small in $S$ if for each non-unit $b$ in $S$, $a+b$ is also non-unit \cite[p. 77]{Golan1999(b)}. As a final phase of \S\ref{sec:idempotentunitssmall}, in Theorem \ref{smallelementsmonoidsemiringsthm}, we show that if $S$ is a semiring and $M$ a positive monoid, then $f=\sum_{i=0}^{n} f_i X^{m_i}$ is small if and only if $f_0$ is small, where $f_i \in S$ and $m_i \in M$ with $m_i < m_{i+1}$ and $m_0 = 0$.

In \S\ref{sec:irreducibleelements}, we discuss some irreducible elements of monoid semirings. The main result of \S\ref{sec:irreducibleelements} is the semiring version of Eisenstein's Irreducibility Criterion (see Theorem \ref{Eisensteincriterion}). Then, we deduce that if $S$ is a weak Gaussian factorial semidomain and $p$ a prime element of $S$ and \[f= a_n X^n + a_{n-1} X^{n-1} + \cdots +a_0\] is a polynomial in $S[X]$ such that the following condition hold:

\begin{enumerate}
	\item $p \nmid a_n$,
	\item $p \mid a_i$, for all $i <n$,
	\item $p^2 \nmid a_0$. 
\end{enumerate}

Then, $f$ cannot be factored into non-constant polynomials in $S[X]$.

\section{Amount semialgebras}\label{sec:amountsemialgebras}

We begin this section by introducing the amount functions.

\begin{definition}[Amount functions] \label{amonutfunctionsdef} Let $S$ be a semiring and $B$ an $S$-semialgebra. We say a function $A$ from $B$ to the set of ideals $\Id(S)$ of $S$, defined by $f \mapsto A_f$, is an amount function if the following properties hold for all $s\in S$ and $f,g\in B$:
	
	\begin{enumerate}
		\item $A$ preserves 0 and 1, i.e. $A_0 = (0)$ and $A_1 = S$.
		\item If $A_f =(0)$, then $f=0$. 
		\item $A$ is homogeneous, i.e. $A_{sf} = sA_f$.
		\item $A$ is submultiplicative, i.e. $A_{fg} \subseteq A_f A_g$.
	\end{enumerate}
	
\end{definition}

\begin{remark}
	Let $S$ be a semiring, $B$ an $S$-semialgebra and $A$ an amount function from $B$ to $\Id(S)$. The homogenity of $A$ implies that \[A_s = (s) \qquad\forall~s\in S.\] In particular, in Definition \ref{amonutfunctionsdef}, the condition $A_0 = (0)$ is superfluous.
\end{remark}

\begin{examples} \label{examplesamountfunctions} In the following, we give two important examples of amount functions:
	
	\begin{enumerate}
		\item Let $(\Gamma,+,0,<)$ be a totally ordered additive monoid and $S$ be a semiring. Let $f=s_1 X^{\alpha_1} + s_2 X^{\alpha_2} + \cdots + s_n X^{\alpha_n}$ be an element of the monoid semiring $S[\Gamma]$. Define the content of $f$, denoted by $c(f)$, to be an ideal of $S$ generated by the coefficients of $f$, i.e. \[c(f) := (s_1,s_2,\ldots,s_n).\] It is easy to verify that $c: S[\Gamma] \longrightarrow \Id(S)$ is an amount function. Note that by $\Id(S)$, we mean the set of all ideals of the semiring $S$.
		
		\item Let us recall that an element $x$ of a totally ordered monoid $(\Gamma,+, <)$ is finitely decomposable if there are only finitely many pairs $(y_i,z_i)$ of elements of $\Gamma$ such that $x=y_i + z_i$. Now, let $(\Gamma, +, 0, <)$ be a totally ordered additive monoid. Also, assume that $\Gamma$ is positive and that each element of $\Gamma$ is finitely decomposable (for example, let $\Gamma = \bigoplus \mathbb N_0$). Let $S$ be a semiring and $S[[\Gamma]]$ be the set of all functions $f: \Gamma \rightarrow S$. Let $f$ and $g$ be arbitrary elements of $S[[\Gamma]]$ and define their addition and multiplication as follows: \[(f+g)(x) = f(x)+g(x), ~ (fg)(x) = \sum_{y+z=x} f(y)g(z).\] It is straightforward to see that $S[[\Gamma]]$ is an $S$-semialgebra \cite{GilmerGramsParker1975}. For each $f\in S[[\Gamma]]$, define $A_f$ to be an ideal of $S$ generated by all $f(s)$, i.e. coefficients of $f$. It is easy to see that the function $A$ from $S[[\Gamma]]$ to $\Id(S)$ is an amount function. For instance, for an element $f=s_0 + s_1 X + \cdots + s_n X^n + \cdots $ in $S[[X]]$, \[A_f = (s_0, s_1, \ldots, s_n, \ldots).\]
	\end{enumerate}	
\end{examples}

In the following, we give another family of examples of amount functions.

\begin{proposition}\label{contentisamount}
Let $B$ be a content $S$-semialgebra. Then, the content function $c$ from $B$ into $\Id(S)$ (see Definition 30 in \cite{Nasehpour2016}) is an amount function.
\end{proposition}

\begin{proof}
By Proposition 31 in \cite{Nasehpour2016}, \[c(f) = \bigcap \{I : I \in \Id(S),~f \in IB\}.\] This simply implies that $c(0) = (0)$. By definition, $c(1) = S$ (see Definition 30 in \cite{Nasehpour2016}). Therefore, $c$ preserves 0 and 1. For the moment, let $c(f) = (0)$. This implies that $f \in (0)B = (0)$. So, $f = 0$. Any content function is, by definition, homogeneous. Finally, by Proposition 31 in \cite{Nasehpour2016}, $c$ is submultiplicative and the proof is complete.
\end{proof}

\begin{definition} \label{amountsemialgebras}
	Let $S$ be a semiring and $B$ an $S$-semialgebra. We say $B$ is an amount $S$-semialgebra if the following conditions hold:
	
	\begin{enumerate}
		\item There is an amount function $A$ from $B$ to $\Id(S)$ defined by $f \mapsto A_f$ with this property that for all $f,g \in B$, there are non-negative integers $m,n$ such that \[A_f ^m A_g ^n A_{fg} = A_f ^{m+1} A_g ^{n+1}\qquad \text{(The Amount Formula).} \label{Amountformula} \]
		
		\item There is a function $\epsilon$ from $\Id(S)$ to $\Id(B)$ defined by $I \mapsto I^{\epsilon}$ with the following properties:
		
		\begin{enumerate}
			\item $A_f \subseteq I$ if and only if $f\in I^{\epsilon}$, for all $f\in B$ and $I\in \Id(S)$.
			\item $I^{\epsilon} \cap S = I$, for all $I\in \Id(S)$.
		\end{enumerate}
	\end{enumerate}
\end{definition}

Amount semialgebras are a generalization of content semialgebras as we show in the following: 

\begin{theorem}
	\label{contentisamount2}
	Let $B$ be a content $S$-semialgebra. Then, $B$ is an amount $S$-semialgebra.
\end{theorem}	

	\begin{proof} Assume that $B$ is a content $S$-semialgebra. By Proposition \ref{contentisamount}, $c$ is an amount function. Obviously, the Dedekind-Mertens content formula (see Definition 30 in \cite{Nasehpour2016}) is a kind of the Amount Formula given in Definition \ref{amountsemialgebras}. Now, define $I^{\epsilon} = IB$. Clearly, $c(f) \subseteq I$ if and only if $f\in IB$ for all $f\in B$ and $I\in \Id(S)$, since $c(f)$ is the smallest ideal satisfying the condition $f\in IB$ (see Definition 24 in \cite{Nasehpour2016}). Finally, it is clear that $I \subseteq IB \cap S$. Now, let $s\in IB \cap S$. So, $c(s) \subseteq I$. However, $c(s) = (s)$ for all $s\in S$. Therefore, $s\in I$. Hence, $ IB \cap S \subseteq I$. From all we said, we conclude that $B$ is an amount $S$-semialgebra and the proof is complete. \end{proof}

\begin{examples}\label{ExamplesAmountSemialgebras}
	Let $S$ be a semiring.
	
	\begin{enumerate}
		\item By Theorem 3 in \cite{Nasehpour2016}, $S[X]$ is an amount (in fact, a content) $S$-semialgebra if and only if $S$ is a subtractive semiring. Note that a semiring $S$ is subtractive if each ideal of $S$ is subtractive. 
		
		\item If $S$ is a Noetherian commutative ring, then by Theorem 2.6 in \cite{EpsteinShapiro2016}, $S[[X]]$ is an amount (in fact, a content) $S$-algebra. For a Dedekind–Mertens lemma for power series in an arbitrary set of indeterminates, check Theorem 3.10 in \cite{GiauToanVo2022}.
		
		\item Let $S$ be a commutative ring and $M$ a (commutative) monoid. Then, by Theorem 2 in \cite{Nasehpour2011}, $S[M]$ is an amount (in fact, a content) $S$-algebra if and only if $M$ is a cancellative and torsion-free monoid. 
		
		\item If $S$ is a valuation ring, then by Theorem 2.8 in \cite{ParkKangToan2018}, for all $f,g \in S[[X]]$, we have \[A^2_f A_g = A_f A_{fg} \text{~or~} A^2_g A_f = A_g A_{gf}\] and in such a case, $S[[X]]$ is an amount $S$-algebra.
		
		\item If $S$ is a Pr\"{u}fer domain, then by Corollary 2.9 in \cite{ParkKangToan2018}, for all $f,g \in S[[X]]$, we have \[(A_f A_g)^2 = A_f A_g A_{fg}\] and again, $S[[X]]$ is an amount $S$-algebra.
	\end{enumerate}
\end{examples}

\begin{remark}
By Theorem \ref{contentisamount2}, we know that any content semialgebra is an amount semialgebra, and because of recent results in commutative algebra, there are some amount (semi)algebras that are not content (semi)algebras (see the last two examples given in Examples \ref{ExamplesAmountSemialgebras}). Recall that the amount (in fact, the content) of any element in a content algebra is finitely generated while this is not necessarily the case in arbitrary amount algebras. Finally, we add that if $I$ is an ideal of a domain $D$ and $D[[X]]$ is an amount $D$-algebra, then by passing from $D$ to the quotient ring $D/I$, the similar amount formulas hold for the amount $D/I$-algebra $(D/I)[[X]]$, and for example, by the last two examples given in Examples \ref{ExamplesAmountSemialgebras}, we obtain new amount algebras with proper zero-divisors. By definition, a zero-divisor is proper if it is nonzero.   
\end{remark}

\begin{definition}
We say an element $f$ of an amount $S$-semialgebra $B$ is primitive if $A_f = S$.
\end{definition}

\begin{proposition}
	Let $B$ be an amount $S$-semialgebra. Then, the following statements hold:
	
	\begin{enumerate}
		\item $fg$ is primitive if and only if $f$ and $g$ are primitive, for all $f$ and $g$ in $B$.
		\item $f\in A_f^{\epsilon}$ for all $f\in B$.
		\item $I \subseteq J$ if and only if $I^{\epsilon} \subseteq J^{\epsilon}$ for all ideals $I$ and $J$ of $S$.
	\end{enumerate}
	
	\begin{proof} 
		(1) Since $A$ is submultiplicative $A_{fg} = S$ implies that $A_f = A_g = S$. On the other hand, if $A_f = A_g = S$, then by the Amount Formula in Definition \ref{amountsemialgebras}, we have $A_{fg} = S$.
		
		(2): Since $A_f \subseteq A_f$, by definition, $f\in A_f^{\epsilon}$.
		
		(3): Assume that $I \subseteq J$ and let $f\in I^{\epsilon}$. By definition, $A_f \subseteq I$. So, $A_f \subseteq J$. This implies that $f\in J^{\epsilon}$. On the other hand, if $I^{\epsilon} \subseteq J^{\epsilon}$, then $I^{\epsilon} \cap S \subseteq J^{\epsilon} \cap S$ which is equivalent to say that $I \subseteq J$.\end{proof}
	
\end{proposition}

Inspired by the definition of reduction in commutative algebra (see Definition 1 in \cite{NorthcottRees1954}), we give the following definition:

\begin{definition}
	Let $I$ and $J$ be ideals of a semiring $S$. We say $J$ is a reduction of $I$ if $J\subseteq I$ and $JI^k = I^{k+1}$ for some positive integer $k$.
\end{definition}

\begin{lemma}
	
	\label{amountreduction}
	Let $B$ be an amount $S$-semialgebra. Then, $A_{fg}$ is a reduction of $A_f A_g$ for all $f,g\in B$.
	
	\begin{proof}  Let $f,g\in B$. Then by definition, there are non-negative integers $m,n$ such that $A_f ^m A_g ^n A_{fg} = A_f ^{m+1} A_g ^{n+1}.$ Let $k=1+\max\{m,n\}$. So, $A_{fg} (A_f A_g) ^{k} = (A_f A_g) ^{k+1}$. Clearly, $k$ is a positive integer and $A_{fg} \subseteq A_f A_g$. Hence, $A_{fg}$ is a reduction of $A_f A_g$ and the proof is complete.\end{proof}
\end{lemma}

\begin{theorem}\label{amountisweakcontent}
	Let $B$ be an amount $S$-semialgebra. Then, $A_f A_g \subseteq \sqrt{A_{fg}}$ for all $f,g\in B$.
	
	\begin{proof} Let $\mathfrak{p}$ be a prime ideal of $S$ containing $A_{fg}$. By Lemma \ref{amountreduction}, $A_{fg}$ is a reduction of $A_f A_g$. So, $A_{fg} (A_f A_g) ^{k} = (A_f A_g) ^{k+1}$ for some positive integer $k$. This implies that $\mathfrak{p}$ contains $A_f A_g$. Hence, in view of Theorem 3.2 in \cite{NasehpourP}, \[\displaystyle A_f A_g \subseteq \bigcap_{\mathfrak{p} \supseteq A_{fg}} \mathfrak{p} = \sqrt{A_{fg}}.\] This completes the proof. \end{proof}
\end{theorem}

\begin{theorem}
	
	\label{primesamountextension}
	
	Let $B$ be an amount $S$-semialgebra. Then, $\mathfrak{p}$ is a prime ideal of $S$ if and only if $\mathfrak{p}^{\epsilon}$ is a prime ideal of $B$.
	
\end{theorem}

\begin{proof}
The proof is the mimicking of the proof of its counterpart in ring theory (see Theorem 15 in \cite{Nasehpour2022}), and so, omitted.
\end{proof}

\section{Semialgebras having Property (A)}\label{sec:PropertyAsemialgebras}

Let us recall that a commutative ring $R$ has Property (A) if any finitely generated ideal $I$ of the ring $R$ with $I \subseteq Z(R)$ is annihilated by a nonzero element $r$ of $R$. It is clear that for any finitely generated ideal $I$ of a ring $R$, one can find a polynomial $f \in R[X]$ such that $I = c(f)$. Based on this point, we give the following definition: 

\begin{definition}\label{PropertyAsemialgebradef}
Let $B$ be an $S$-semialgebra and $A: B \rightarrow \Id(S)$ an amount function. We say the $S$-semialgebra $B$ has Property (A) if for all $f \in B$, the inclusion $A_f \subseteq Z(S)$ implies that there is a nonzero $s \in S$ with $s A_f = (0)$.
\end{definition}

Inspired by the definition of the ``prime avoidance'' property in commutative ring theory, given by Chen in \cite{Chen2021}, we give the following definition:

\begin{definition}\label{primeavoidancedef}
	Let $S$ be a semiring and $\Gamma \subseteq \Spec(S)$. We say $\Gamma$ satisfies the prime avoidance property if for any ideal $I$ of $S$ the following condition holds:
	
	\begin{itemize}
		\item $I \subseteq \bigcup_{\mathfrak{p} \in \Gamma} \mathfrak{p}$ implies that $I \subseteq \mathfrak{p}$ for some $\mathfrak{p} \in \Gamma$.
	\end{itemize}
\end{definition}

\begin{examples}\label{PrimeAvoidanceExamples} In the following, we give some examples for the prime avoidance property for rings and semirings:
	
	\begin{enumerate}
		\item Let $S$ be a semiring and $\Gamma \subseteq \Spec(S)$. If $\Gamma$ is a finite set and each element of $\Gamma$ is a subtractive prime ideal of $S$, then $\Gamma$ satisfies the prime avoidance property and its proof is similar to the proof of the Prime Avoidance Theorem given in Theorem 81 in \cite{Kaplansky1974}.
		
		\item Let $R$ be a ring such that each prime ideal of $R$ is the radical of a principal ideal in $R$. Then, by the main theorem of the paper \cite{Smith1971}, each $\Gamma \subseteq \Spec(R)$ satisfies the prime avoidance property. 
	\end{enumerate}	
\end{examples}

\begin{definition}\label{controllingzerodivisors}
Let $S$ be a semiring and $\Gamma \subseteq \Spec(S)$. We say $Z(S)$ is controlled by $\Gamma$ if $Z(S) = \bigcup_{\mathfrak{p} \in \Gamma} \mathfrak{p}$.
\end{definition}

\begin{remark}
\begin{enumerate}
	\item By definition, a commutative ring $R$ has few zero-divisors if $Z(R)$ is a finite union of prime ideals \cite{Huckaba1988}. So, $R$ has few zero-divisors if and only if $Z(R)$ is controlled by a finite number of prime ideals.
	
	\item It is evident that a semiring $S$ has very few zero-divisors if and only if $Z(S)$ is controlled by a finite number of prime ideals in $\Ass(S)$. Let us recall that a semiring $S$ has very few zero-divisors if the set of zero-divisors of $S$ is a finite union of prime ideals in $\Ass(S)$ (see Definition 48 in \cite{Nasehpour2016}).
\end{enumerate}

\end{remark}

\begin{theorem}\label{PropertyAsemialgebrathm}
Let $S$ be a semiring and $\Gamma \subseteq \Ass(S)$. If $Z(S)$ is controlled by $\Gamma$ and $\Gamma$ satisfies the prime avoidance property, then any amount $S$-semialgebra has Property (A). 
\end{theorem}

\begin{proof}
Let $B$ be an amount $S$-semialgebra and $f$ be an element of $B$ with $A_f \subseteq Z(S)$. By assumption \[Z(S) = \bigcup_{\mathfrak{p} \in \Gamma} \mathfrak{p}\] and $\Gamma \subseteq \Ass(S)$ satisfies the prime avoidance property. This implies that $A_f \subseteq \mathfrak{p}$, where $\mathfrak{p} = \Ann(s)$, for some $s\in S$. So, $s A_f = (0)$. Hence, the $S$-semialgebra $B$ has Property (A), as required.
\end{proof}

In view of Examples \ref{PrimeAvoidanceExamples}, we have the following:

\begin{corollary}\label{PropertyAsemialgebracor}
If $S$ has very few zero-divisors (for example, $S$ is Noetherian), then the $S$-semialgebra $S[X]$ has Property (A).
\end{corollary}

The definition of graded semiring is similar to its counterpart in ring theory. One can check the paper \cite{AllenKimNeggers2013} for some details of graded semirings. In this paper, by a graded semiring, we mean a semiring graded by the integers. A graded semiring is nontrivial in case it contains a regular homogeneous element of nonzero degree.

\begin{theorem}
Let $S$ be a nontrivial graded semiring and $X$ an indeterminate over $S$. Then, the $S$-semialgebra $S[X]$ has Property (A).
\end{theorem}

\begin{proof}
In view of McCoy's property on arbitrary semirings proved in Corollary 1.2 in \cite{Nasehpour2021}, the proof is similar to the proof of Theorem 2.7 in \cite{Huckaba1988}, and so, omitted. 
\end{proof}

\section{Armendariz and McCoy semialgebras}\label{sec:armendarizmccoysemialgebras}

Let $R$ be a ring. Armendariz, in Lemma 1 in \cite{Armendariz1974}, proves that if $R$ is a reduced ring and $f$ and $g$ are elements of the polynomial ring $R[X]$ with $f = \sum_{i=1}^{m} a_i X^i$ and $g = \sum_{j=1}^{n} b_j X^j$, then $fg = 0$ implies $a_i b_j = 0$, for all $i$ and $j$. A ring $R$ is an Armendariz ring if whenever polynomials $f = \sum_{i=1}^{m} a_i X^i$ and $g = \sum_{j=1}^{n} b_j X^j$ in $R[X]$ satisfy $fg = 0$, we have $a_i b_j = 0$ for every $i$ and $j$ (see Definition 1.1 in \cite{RegeChhawchharia1997}). In other words, if $R$ is a commutative ring, then $R$ is Armendariz if and only if $fg = 0$ implies that $c(f) c(g) = (0)$ for all $f$ and $g$ in $R[X]$. Now, let $B$ be an $R$-algebra such that $B$ is a content $R$-module. It is said that $B$ is an Armendariz $R$-algebra if $fg = 0$ implies $c(f) c(g) = (0)$, for all $f$ and $g$ in $B$ (check Definition 6 in \cite{Nasehpour2016AB}). 

Similarly, we define Armendariz semialgebras. Let $S$ be a semiring, $B$ an $S$-semialgebra, and $A$ an amount function from $B$ into $\Id(S)$. Using Definition \ref{amonutfunctionsdef}, it is easy to see that the equality $A_f A_g = (0)$ implies $fg=0$ for all $f,g \in B$. Based on this fact, we give the following definition:

\begin{definition}\label{Armendarizsemialgebradef}
	Let $B$ be an $S$-semialgebra and $A: B \rightarrow \Id(S)$ an amount function. We say $B$ is an Armendariz $S$-semialgebra if for all $f,g \in B$, $fg = 0$ implies that $A_f A_g = (0)$.	
\end{definition}

\begin{proposition} \label{AmountisArmendariz}
	Let $S$ be a nilpotent-free semiring and $B$ an amount $S$-semialgebra. Then, $B$ is Armendariz.
\end{proposition}	
	\begin{proof}
		Let $f$ and $g$ be elements of $B$ such that $fg=0$. By the Amount Formula in Definition \ref{amountsemialgebras}, there are non-negative integers $m$ and $n$ such that \[A_f ^{m+1} A_g ^{n+1}=(0). \] Since $S$ is nilpotent-free, $A_f A_g = (0)$. Hence, $B$ is Armendariz, as required.
	\end{proof}

The following is a semiring version of Proposition 3.5 in \cite{GilmerGramsParker1975}:

\begin{proposition} \label{proposition35GGP}
Let $S$ be a semiring such that each minimal prime ideal of $S$ is subtractive. Also, let $f = \sum_{i=0}^{+\infty} f_i X^i$ and $g = \sum_{j=0}^{+\infty} g_j X^j$ be elements of $S[[X]]$ such that $fg = 0$. Then, $A_f A_g \subseteq \Nil(S)$.
\end{proposition}

\begin{proof}
On the contrary, assume that $A_f A_g$ is not a subset of $\Nil(S)$. Since $\Nil(S)$ is the intersection of all minimal prime ideal of $S$, we see that $A_f A_g$ is not a subset of $\mathfrak{p}$ for some minimal prime ideal $\mathfrak{p}$ of $S$. This implies that both $A_f$ and $A_g$ are not subset of $\mathfrak{p}$. Let $m$ and $n$ be smallest indices of $f$ and $g$, respectively, that $f_m$ and $g_n$ are not elements of $\mathfrak{p}$. Assume that $h_k$s are the coefficients of $fg$, i.e.,  \[fg = h = \sum_{k=0}^{+\infty} h_k X^k.\] It is clear that \[h_{m+n} = f_0 g_{m+n} + \dots + f_m g_n + \dots + f_{m+n} g_0.\] Since $fg = 0$, $h_{m+n} = 0 \in \mathfrak{p}$. This implies that $f_m g_n \in \mathfrak{p}$ because $f_i$ and $g_j$ are elements of $\mathfrak{p}$ if $i < m$ and $j < n$ and $\mathfrak{p}$ is subtractive. However, this is a contradiction since $f_m$ and $g_n$ are not elements of the prime ideal $\mathfrak{p}$. Hence, $A_f A_g \subseteq \Nil(S)$, as required.
\end{proof}

Inspired by Proposition \ref{proposition35GGP}, we give the following definition:

\begin{definition}\label{GGPsemialgebradef} Let $B$ be an $S$-semialgebra and $A$ an amount function from $B$ into $\Id(S)$. We say that $B$ is a Gilmer Grams Parker (for short, GGP) $S$-semialgebra if $fg = 0$ implies that $A_f A_g \subseteq \Nil(S)$, for all $f$ and $g$ in $B$.
\end{definition}

\begin{corollary}
If $B$ is a GGP $S$-semialgebra and $S$ is nilpotent-free, then $B$ is an Armendariz $S$-semialgebra. In particular, if $S$ is nilpotent-free and each minimal prime ideal of $S$ is subtractive, then $S[[X]]$ is an Armendariz $S$-semialgebra. 
\end{corollary}

McCoy proved that if $R$ is a ring and $fg = 0$ for some $f$ and $g$ in $R[X]$ with $g \neq 0$, then there is a nonzero $r\in R$ such that $fr = 0$ \cite{McCoy1942}. Note that the equality $fr = 0$ is equivalent to $r c(f) = (0)$. Now, we give the following definition which is a generalization of the definition of McCoy semialgebras given in Definition 3.4 in \cite{Nasehpour2021}:

\begin{definition}\label{McCoysemialgebradef}
	Let $B$ be an $S$-semialgebra and $A: B \rightarrow \Id(S)$ an amount function. We say $B$ is a McCoy $S$-semialgebra if for all $f,g \in B$, $fg = 0$ together with $g \neq 0$ implies that there is a nonzero $s \in S$ such that $s A_f = (0)$.
\end{definition}

\begin{proposition}\label{armendarizismccoy}
Any Armendariz $S$-semialgebra is a McCoy $S$-semialgebra.
\end{proposition}

\begin{proof}
Let $fg = 0$ with $g \neq 0$. So, we have $A_f A_g = (0)$ because $B$ is Armendariz. Choose a nonzero element $s$ in $A_g$. Then, $s A_f = (0)$ and the proof is complete. 
\end{proof}

\begin{examples} Let $S$ be a semiring.
	
	\begin{enumerate}
		
		\item By Corollary 2.2 in \cite{Nasehpour2021} if $G$ is a cancellative torsion-free commutative monoid, then the $S$-semialgebra $S[G]$ is a McCoy $S$-semialgebra.
		
		\item By Proposition 31 in \cite{Nasehpour2016}, any content $S$-semialgebra is McCoy but if $S$ is not subtractive, then $S[X]$ is a McCoy $S$-semialgebra (see Corollary 2.2 in \cite{Nasehpour2016}) while by Theorem 3 in \cite{Nasehpour2016} it is not a content $S$-semialgebra.
		
		\item In Question 2.12 \cite{Nasehpour2018}, the author asked if there was any faithfully flat McCoy algebra that was not a content algebra. Interestingly, Epstein found faithfully flat McCoy algebras that are not content algebras (see Example 1 in \cite{Epstein2022}).
	\end{enumerate}	
	
\end{examples}

\begin{theorem}\label{powerseriesismccoy}
Let $S$ be a semiring. The $S$-semialgebra $S[[X]]$ is a McCoy $S$-semialgebra if one of the following conditions hold:

\begin{enumerate}
	\item The semiring $S$ is nilpotent-free and weak Gaussian.
	
	\item The semiring $S$ is zerosumfree.
\end{enumerate}
\end{theorem}

\begin{proof} 
For the proof of the first case, since $S$ is weak Gaussian, for all $f$ and $g$ in $S[[X]]$, we have $A_f A_g \subseteq \sqrt{A_{fg}}$ \cite[Theorem 45]{Nasehpour2016}. Now, if $fg=0$, then $A_f A_g =(0)$ because $S$ is nilpotent-free. Therefore, if we suppose that $g$ is nonzero, we can take $s$ to be a nonzero element in $A_g$ and $s$ annihilates $f$.

For the proof of the second case, let $f=\sum_{i=0}^{\infty} f_i X^i$ be a zero-divisor on $S[[X]]$. So, there is a nonzero $g=\sum_{j=0}^{\infty} g_j X^j$ in $S[[X]]$ such that $fg=0$. This implies that $\sum_{i+j=n} f_i g_j=0$ for each non-negative integer $n$. Since $S$ is zerosumfree, $f_i g_j=0$ for all $i,j\geq 0$. From this we obtain that $A_f A_g = (0)$. Since $g$ is nonzero, $A_g$ is nonzero. Take $s$ to be a nonzero element in $A_g$. It is obvious that $s$ annihilates $f$ and the proof is complete.\end{proof}

Let us recall that an ideal $I$ of a semiring $S$ is regular if it contains a regular element, i.e. an element that is not a zero-divisor. The following is a generalization of Theorem 13 in \cite{Nasehpour2010}:

\begin{theorem}\label{WhenMcCoyhasPropertyA}
	Let $B$ be a McCoy $S$-semialgebra and $A$ an amount function from $B$ into $\Id(S)$. Then, the following statements are equivalent:
	
	\begin{enumerate}
		\item The $S$-semialgebra $B$ has Property (A).
		\item For all $f \in B$, $f$ is a regular element of $B$ if and only if $A_f$ is a regular ideal of $S$.
	\end{enumerate}

\end{theorem}
	
	\begin{proof}
		$(1) \implies (2)$: Let $B$ have Property (A). If $f \in B$ is regular, then for all nonzero $s \in S$, $sf \not= 0$, and so, for all nonzero $s \in S$, $sA_f \not= (0)$, i.e. $\Ann(A_f) = (0)$ and by the definition of Property (A), $A_f \not\subseteq Z(S)$. This means that $A_f$ is a regular ideal of $S$. Now, let $A_f$ be a regular ideal of $S$. By regularity of $A_f$, we have $A_f \not\subseteq Z(S)$ and so, $\Ann(A_f) = (0)$. This means that for all nonzero $s \in S$, $sA_f \not= (0)$. Consequently, for all nonzero $s \in S$, $sf \not= 0$. Since $B$ is a McCoy $S$-semialgebra, $f$ is not a zero-divisor on $B$, i.e. it is regular.
		
		$(2) \implies (1)$: Choose $f$ in $B$ such that $A_f \subseteq Z(S)$. Since $A_f$ is not a regular ideal of $S$, $f$ is not a regular element of $B$. By assumption, $f$ is a zero-divisor on $B$. Since $B$ is a McCoy $S$-semialgebra, there exists a nonzero $s \in S$ with $sf = 0$. Finally, since $A$ is a homogeneous, we have $sA_f = (0)$ showing that $S$-semialgebra $B$ has Property (A) and the proof is complete.
	\end{proof}

\begin{theorem}\label{moreexamplesarmendarizzemialgebras}
Let $S$ be a zerosumfree semiring. Then, the following statements hold:

\begin{enumerate}
	\item \label{monoidsemiringarmendariz}  If $M$ is a monoid, then $S[M]$ is an Armendariz $S$-semialgebra.
	\item \label{powerseriessemiringarmendariz} If $\{X_i\}_{i=1}^{n}$ is a family of indeterminates over $S$, then the formal power series semiring $S[[X_1,\dots,X_n]]$ is an Armendariz $S$-semialgebra.  
\end{enumerate}
\end{theorem}

\begin{proof}
(1): Let $f=\sum_{i=1}^{k} f_iX^{m_i}$ and $g=\sum_{j=1}^{l} g_jX^{n_j}$ be elements of $S[M]$ such that $fg=0$. Since $S$ is zerosumfree, $S[M]$ is also zerosumfree. This implies that $f_i g_j =0$ for each $1\leq i \leq k$ and $1\leq j \leq l$. Hence, $c(f)c(g) = 0$ and so, $S[M]$ is an Armendariz $S$-semialgebra.

(2): The proof of this statement that $S[[X]]$ is an Armendariz $S$-semialgebra, is similar to the proof of (1) and so, omitted. Finally, an obvious induction gives the desired result and the proof is complete.
\end{proof}

\section{Zero-divisors and nilpotents}\label{sec:zerodivisorsnilpotents}

Let $R$ be a commutative ring and $\{X_i\}_{i\in \Lambda}$ be a set of indeterminates over $R$. In Theorem 9 in \cite{Nasehpour2010}, it is proved that if $Z(R) = \bigcup_{i=1}^{n} \mathfrak{p}_i$, where $\mathfrak{p}_i \in \Ass(R)$, then \[Z(R[\{X_i\}_{i\in \Lambda}]) = \bigcup_{i=1}^{n} \mathfrak{p}_i [\{X_i\}_{i\in \Lambda}].\] We generalize this in the following:

\begin{theorem}\label{extensionofzerodivisors1} 
	Let $S$ be a semiring and $\Gamma \subseteq \Ass(S)$ satisfy the prime avoidance property. Also, let $B$ be an amount $S$-semialgebra, where $A$ is its amount function from $B$ into $\Id(S)$. If $Z(S)$ is controlled by $\Gamma$ and $B$ is a McCoy $S$-semialgebra, then $Z(B)$ is controlled by $\mathfrak{p}^{\epsilon}$s, where $\mathfrak{p} \in \Gamma$.
\end{theorem}

\begin{proof}
	Let $Z(S) = \bigcup_{\mathfrak{p} \in \Gamma} \mathfrak{p}$. We will show that $Z(B) = \bigcup_{\mathfrak{p} \in \Gamma} \mathfrak{p}^{\epsilon}$. Let $f \in Z(B)$. Since $B$ is a McCoy $S$-semialgebra, there is a nonzero $s \in S$ with $sA_f = (0)$. This implies that $A_f \subseteq Z(S)$. On the other hand, $\Gamma \subseteq \Ass(S)$ satisfies the prime avoidance property. Therefore, we have $A_f \subseteq \mathfrak{p}$, for some $\mathfrak{p} \in \Gamma$. Consequently, $f \in \mathfrak{p}^{\epsilon}$. Now, let $f \in \mathfrak{p}^{\epsilon}$, where $\mathfrak{p} \in \Gamma$. This implies that $A_f \subseteq \mathfrak{p}$. Since $\Gamma \subseteq \Ass(S)$, $A_f$ has a nonzero annihilator and this means that $f$ is a zero-divisor of $B$. Hence, $Z(B) = \bigcup_{\mathfrak{p} \in \Gamma} \mathfrak{p}^{\epsilon}$, as required.
\end{proof}

Now, we proceed to prove another result for controlling the zero-divisors of McCoy semialgebras:

\begin{theorem}\label{extensionofzerodivisors2}
	Let $S$ be a semiring and $\Gamma \subseteq \Spec(S)$ satisfy the prime avoidance property. Also, let $B$ be an amount $S$-semialgebra, where $A$ is its amount function from $B$ into $\Id(S)$. If $Z(S)$ is controlled by $\Gamma$ and the $S$-semialgebra $B$ is McCoy and has Property (A), then $Z(B)$ is controlled by $\mathfrak{p}^{\epsilon}$s, where $\mathfrak{p} \in \Gamma$.
\end{theorem}

\begin{proof}
Let $f \in Z(B)$. Since $B$ is a McCoy $S$-semialgebra, there is a nonzero $s\in S$ with $s A_f = (0)$. This implies that $A_f \subseteq Z(S)$. On the other hand, $Z(S) = \bigcup_{\mathfrak{p} \in \Gamma} \mathfrak{p}$ and $\Gamma \subseteq \Ass(S)$ satisfies the prime avoidance property. So, $A_f \subseteq \mathfrak{p}$, for some $\mathfrak{p} \in \Gamma$. From this we obtain that $f \in \mathfrak{p}^{\epsilon}$.

Conversely, let $f \in \mathfrak{p}^{\epsilon}$. So, $A_f \subseteq \mathfrak{p} \subseteq Z(S)$. Since the $S$-semialgebra $B$ has Property (A), the ideal $A_f$ is annihilated by a nonzero $s \in S$. This means that $f$ is a zero-divisor of $B$ and the proof is complete.
\end{proof}

Let us recall that an element $r$ of a ring $R$ is said to be prime to an ideal $I$ of $R$ if $I : (r) =I$, where by $I : (r)$, we mean the set of all elements $a$ of $R$ such that $ar \in I$ (check p. 223 in \cite{ZariskiSamuel1958}). Let $I$ be an ideal of $R$ and $S(I)$ the set of all elements of $R$ that are not prime to $I$. Let us recall that an ideal $I$ of $R$ is said to be primal if $S(I)$ forms an ideal of $R$ \cite{Fuchs1950}. It is easy to see that if $S(I)$ is an ideal of $R$, then it is a prime ideal. A ring $R$ is said to be primal if $(0)$ is a primal ideal. Note that the ring $R$ is primal if and only if $Z(R)$ is an ideal of $R$. Based on this, we give the following definition:

\begin{definition}\label{primalsemiringdef}
We say a semiring $S$ is primal if $Z(S)$ is an ideal of $S$.
\end{definition}

\begin{remark}
It is easy to see that if $S$ is a primal semiring, then $Z(S)$ is a prime ideal of $S$. In other words, a semiring is primal if and only if $Z(S)$ is controlled by a single prime ideal of $S$.
\end{remark}

\begin{corollary}\label{controllingprimalsemirings}
Let $S$ be a semiring and $B$ be an amount $S$-semialgebra, where $A$ is its amount function from $B$ into $\Id(S)$. Also, let the $S$-semialgebra $B$ be McCoy and have Property (A). Then, if the semiring $S$ is primal, then so is the semiring $B$.	
\end{corollary}

\begin{proof}
	If $Z(S) = \mathfrak{p}$, then by Theorem \ref{extensionofzerodivisors2}, $Z(B) = \mathfrak{p}^{\epsilon}$.
\end{proof}

Let us recall that an element $s$ of a semiring $S$ is nilpotent if $s^n = 0$ for some positive integer $n$. We collect all nilpotent elements of $S$ in the set $\Nil(S)$. 

Now, let $M$ be a monoid and $S$ a semiring. One interesting question is to determine all nilpotent elements of $S[M]$. In this direction, we give the following notation: \[\Nil(S)[M] := \left\{\sum_{i=1}^{n} f_i X^{m_i}: f_i \in \Nil(S),~m_i \in M,~n\in \mathbb N\right\}.\]

\begin{theorem}\label{nilpotentsmonoidsemirings}
	Let $M$ be a monoid and $S$ a semiring. Then, the following statements hold:
	
	\begin{enumerate}
		\item Each element of $\Nil(S)[M]$ is nilpotent.
		
		\item $\Nil(S[M]) = \Nil(S)[M]$ if one of the following conditions holds:
		
		\begin{enumerate}
			\item $S$ is a zerosumfree semiring.
			\item $S$ is a weak Gaussian semiring and $M$ is a totally ordered monoid.
		\end{enumerate}
	\end{enumerate}
	
	\begin{proof} (1) The proof is straightforward. 
		
		(2)(a): Let $S$ be zerosumfree and $f$ a nonzero and nilpotent element of $S[M]$. Write $f = \sum_{i=1}^{n} f_i X^{m_i}$ and suppose that there is a positive integer $k$ such that $f^k = 0$. Since $S$ is zerosumfree, each summand appearing in each coefficient of $f^k$ is zero. In particular, $f_{i}^{k} =0$ for each $i$. Consequently, $f\in \Nil(S)[M]$.
		
		(2)(b): Let there be a positive integer $k$ such that $f^k =0$. Let $\mathfrak{p}$ be an arbitrary prime ideal of $S$. Obviously, $f^k\in \mathfrak{p}[M]$. Since $S$ is weak Gaussian, each prime ideal $\mathfrak{p}$ of $S$ is subtractive. Therefore, by Theorem 2.3 in \cite{Nasehpour2021}, $\mathfrak{p}[M]$ is prime, and so, $f\in \mathfrak{p}[M]$ which is equivalent to say that $c(f) \subseteq \mathfrak{p}$. Consequently, $c(f) \subset \bigcap_{\mathfrak{p}\in \Spec(S)} \mathfrak{p}$. By Proposition 7.28 in \cite{Golan1999(b)}, $\bigcap_{\mathfrak{p}\in \Spec(S)} \mathfrak{p} = \Nil(S)$. Hence, $f\in \Nil(S)[M]$ and the proof is complete. \end{proof}
\end{theorem}

A semiring $S$ is called nilpotent-free if $s^2 = 0$ implies $s = 0$, for all $s\in S$.

\begin{corollary}
	Let $S$ be a zerosumfree semiring and $M$ a monoid. Then, $S$ is nilpotent-free if and only if $S[M]$ is so.
\end{corollary}

It is clear that $\Nil(S) \subseteq Z(S)$. Based on this, we give the following definition:

\begin{definition}\label{entirelikesemiringsdef}
We say a semiring $S$ is entire-like if $Z(S) \subseteq \Nil(S)$.
\end{definition}

\begin{proposition}\label{monoidsemiringentirelikepro}
	Let $M$ be a monoid and $S$ an entire-like semiring. If $S[M]$ is a McCoy $S$-semialgebra, then $S[M]$ is entire-like.
	
	\begin{proof} Let $f$ be a zero-divisor on $S[M]$. Since $S[M]$ is a McCoy $S$-semialgebra, there is a nonzero element $c\in S$ such that $c$ annihilates each coefficient of $f$. This means that each coefficient of $f$ is a zero-divisor on $S$. Since $S$ is entire-like, each coefficient of $f$ is nilpotent. This implies that $f$ is nilpotent. Hence, $S[M]$ is entire-like.
	\end{proof}
\end{proposition}

Now, we proceed to determine conditions under which $S[M]$ is free of nontrivial zero-divisors, i.e., it is entire. Recall that a semiring $S$ is an information algebra if it is an entire zerosumfree semiring \cite[p. 4]{Golan1999(b)}.

\begin{theorem}
	Let $M$ be a monoid and $S$ a semiring. If $S$ is an information algebra, then so is $S[M]$.
	
	\begin{proof} Let $f$ and $g$ be nonzero elements of $S[M]$. Set 	$f=\sum_{i=1}^{m} f_i X^{s_i}$ and $g=\sum_{j=1}^{n} g_j X^{t_j}$, where all $f_i$s and $g_j$s are nonzero. Since $S$ is entire, $f_i g_j$ for all $1\leq i \leq m$ and $1\leq j \leq n$ are nonzero. Now, since $S$ is zerosumfree, $fg$ cannot be zero. Hence, $S[M]$ is an information algebra, as required.\end{proof}
\end{theorem}

Since any bounded distributive lattice is an example of an information algebra, we have the following:

\begin{corollary}
	Let $M$ be a monoid and $L$ a bounded distributive lattice. Then, $L[M]$ is an information algebra.
\end{corollary}

Let $\Id(S)$ be the set of all ideals of a semiring $S$. Define addition and multiplication of ideals of $S$ as follows: \[I+J=\{x+y: x\in I,~y\in I\},\] and \[IJ=\left\{\sum_{i=1}^{n} x_i y_i: x_i\in I,~y_i\in J,~n\in \mathbb N\right\}.\] Then, $\Id(S)$ with the above addition and multiplication becomes a zerosumfree semiring \cite[Proposition 6.29]{Golan1999(b)}.

\begin{corollary} Let $M$ be a monoid and $S$ a semiring. Then, the following statements are equivalent:

\begin{enumerate}
	\item The semiring $S$ is entire.
	\item The semiring $\Id(S)$ is an information algebra.
	\item The semiring $\Id(S)[M]$ is an information algebra.
\end{enumerate}
\end{corollary}

\begin{definition}
	
	\label{degreemonoidsemiringdef}
	
	Let $M$ be a totally ordered monoid and set $M_{\infty} = M \cup \{-\infty\}$. Define $-\infty < m$ for all $m\in M$ and $-\infty + m = m + (-\infty) = -\infty$ for all $m\in M_{\infty}$. If $S$ is a semiring, the order and degree of an element $f$ in $S[M]$, denoted by $\ord(f)$ and $\deg(f)$, respectively, are defined as follows:
	
	\begin{enumerate}
		\item If $f=0$, then $\deg(f) = -\infty$,
		\item If $f$ is a monomial, i.e. $f=aX^m$, where $a\neq 0$ is in $S$ and $m\in M$, then $\ord(f) = \deg(f)=m$. 
		\item If $f = \sum_{i=1}^{n} a_i X^{m_i}$, where $a_i$s are nonzero elements of $S$ and $m_i$s are in $M$, $n\geq 2$, and $m_i < m_{i+1}$ for all $1\leq i\leq n-1$, then $\ord(f) = m_1$ and $\deg(f)=m_n$. 
	\end{enumerate} 
\end{definition}

The proof of the following is straightforward:

\begin{proposition}
	
	\label{degreemonoidsemiringpro}
	
	Let $S$ be a semiring and $M$ a totally ordered monoid. Then, for all $f$ and $g$ in $S[M]$, we have:
	
	\begin{enumerate}
		\item $\deg(f+g) \leq \max \{\deg(f),\deg(g)\}$.
		\item $\deg(fg) \leq \deg(f) + \deg(g)$.
		\item If $f$, $g$, and $f+g$ are nonzero, then $\ord(f+g) \leq \min \{\ord(f),\ord(g)\}$.
		\item If $f$, $g$, and $fg$ are nonzero, then $\ord(f)+\ord(g) \leq \ord(fg)$.
	\end{enumerate}
	
\end{proposition}

\begin{theorem}
	
	\label{monoidsemiringzerosumfree}
	
	Let $S$ be a semiring and $M$ a totally ordered monoid. Then, the following statements are equivalent:
	
	\begin{enumerate}
		\item For all $f$ and $g$ in $S[M]$, we have $\deg(f+g) = \max \{\deg(f),\deg(g)\}$.
		\item The semiring $S$ is zerosumfree.
		\item The semiring $S[M]$ is zerosumfree.
	\end{enumerate}
	
	\begin{proof}
		
		$(1) \rightarrow (2)$: Let $a$ and $b$ be in $S$ such that $a+b=0$. So, $\deg((a+b)X^m) = -\infty$. This implies that $\deg(aX^m) = \deg(bX^m) = -\infty$ and this means that $a=b=0$. So, $S$ is zerosumfree. 
		
		$(2) \rightarrow (3)$ is obvious.
		
		$(3) \rightarrow (1)$: Let $f$ and $g$ be in $S[M]$. If both $f$ and $g$ are nonzero and $\deg(f) = \deg(g)$, then $f+g$ cannot be zero because $S$ is zerosumfree. So, \[\deg(f+g) = \max \{\deg(f),\deg(g)\}.\] The other cases are obvious and the proof is complete. \end{proof}
	
\end{theorem}

\begin{example}
	Since for any semiring $S$, the semiring $\Id(S)$ is zerosumfree \cite[Proposition 6.29]{Golan1999(b)}, this gives plenty of examples for Theorem \ref{monoidsemiringzerosumfree}.
\end{example}

\begin{theorem}
	
	\label{monoidsemiringentire}
	
	Let $S$ be a semiring and $M$ a totally ordered monoid. Then, the following statements are equivalent:
	
	\begin{enumerate}
		\item For all $f$ and $g$ in $S[M]$, we have $\deg(fg) = \deg(f)+\deg(g)$.
		\item The semiring $S$ is entire.
		\item The semiring $S[M]$ is entire.
	\end{enumerate}
	
	\begin{proof}
		
		$(3) \rightarrow (1)$: If either $f$ or $g$ is zero, then the statement holds obviously. Now, let $f=\sum_{i=1}^{k} s_i X^{m_i}$ and $g=\sum_{j=1}^{l} t_j X^{n_j}$ such that $m_i < m_{i+1}$ and $n_j < n_{j+1}$ for all $i$ and $j$, and $s_k$ and $t_l$ are nonzero in $S$. It is clear that $s_k t_l$ is nonzero and the monomial $s_k t_l X^{m_k + n_l}$ appearing in $fg$ has the highest degree.
		
		$(1) \rightarrow (2)$: Let $a$ and $b$ be both nonzero elements of $S$. It is clear that $\deg(abX^0) = \deg(aX^0) \deg(bX^0) = 0$. So, $ab$ is not zero and $S$ is entire.
		
		$(2) \rightarrow (3)$: The proof of this implication is similar to the proof of Theorem 8.1 in \cite{Gilmer1984} and therefore, omitted.\end{proof}
	
\end{theorem}

By Theorem \ref{monoidsemiringzerosumfree} and Theorem \ref{monoidsemiringentire}, we have the following result:

\begin{corollary}
	
	\label{monoidsemiringzerosumfreeentire}
	
	Let $S$ be a semiring and $M$ a totally ordered monoid. Then, the following statements are equivalent:
	
	\begin{enumerate}
		\item The function $\deg: (S[M],+,\cdot,0,1) \longrightarrow (M_{\infty},\max,+,-\infty,0)$ is a semiring morphism.
		\item The semiring $S$ is an information algebra.
		\item The semiring $S[M]$ is an information algebra.
	\end{enumerate} 	 
\end{corollary}

\begin{theorem}
	
	\label{ordsemiring}
	
	Let $S$ be an entire semiring and $M$ a totally ordered monoid. Then, the following statements hold for all $f$ and $g$ in $S[M]$:
	
	\begin{enumerate}
		\item If $f$ and $g$ are nonzero, then $\ord(fg) = \ord(f) + \ord(g)$.
		\item If $fg$ is a nonzero monomial element of $S[M]$, then $f$ and $g$ are both monomials.
	\end{enumerate}

	\begin{proof} The proof of Statement (1) is similar to the proof of the implication $(3) \rightarrow (1)$ in Theorem \ref{monoidsemiringentire} and so, omitted.
		
		(2):  Let $fg\neq 0$ be a monomial. Then, \[\ord(f) + \ord(g) = \ord(fg) = \deg(fg) = \deg(f)+\deg(g).\] This implies that $\ord(f) = \deg(f)$ and $\ord(g) = \deg(g)$. Hence, $f$ and $g$ are both monomials and the proof is complete.
	\end{proof}
\end{theorem}

\section{$k$-potents, units, and small elements}\label{sec:idempotentunitssmall}

Let us recall that an element $x$ in a semiring $S$ is multiplicatively idempotent if $x^2 = x$ \cite[p. 3]{Golan1999(b)}. Let $k\geq 2$ be a positive integer. An element $x$ in a ring $R$ is $k$-potent if $x^k = x$ \cite{Mosic2015}. Based on these, we define an element $x$ in a semiring $S$ to be multiplicatively $k$-potent if $x^k =x$. An element $m$ of an additive monoid $M$ is $k$-potent if $k\cdot m = m$.

\begin{proposition}\label{exampleskpotent} Let $S$ be a semiring and $M$ a monoid. Then, the following statements hold:
	
	\begin{enumerate}
		
		\item If $s$ is a multiplicatively $k$-potent element of $S$ and $m$ is a $k$-potent element of $M$, then $sX^m$ is a multiplicatively $k$-potent element of $S[M]$.
		
		\item If each element of the set $\{s_i\}_{i=1}^n \subseteq S$ is multiplicatively $k$-potent such that $s_i s_j = 0$ for all $i\neq j$, then $\sum_{i=1}^{n} s_i$ is also multiplicatively $k$-potent.
		
		\item If for each $1\leq i \leq n$, $s_i$ in $S$ is multiplicatively $k$-potent, $s_i s_j = 0$ for all $1 \leq i < j \leq n$, and for each $1\leq i \leq n$, $m_i$ in $M$ is $k$-potent, then $\sum_{i=1}^{n} s_i X^{m_i}$ is a multiplicatively $k$-potent element of $S[M]$.
		
		\item If $\{m_i\}_{i=1}^n$ is a subgroup of the monoid $M$ and $n\cdot 1=\sum_{i=1}^{n} 1$ is a unit in $S$, then $\displaystyle \frac{\sum_{i=1}^{n} X^{m_i}}{n\cdot 1}$ is a multiplicatively idempotent element of $S[M]$.  
	\end{enumerate}
\end{proposition}	
\begin{proof} The statement (1) is straightforward. In view of the Multinomial Theorem, the proof of the statements (2) and (3) is also straightforward. 
	
(4): Let $f= \sum_{i=1}^{n} X^{m_i}$. Since $\{m_i\}_{i=1}^n$ is a subgroup of the monoid $M$, $X^{m_j} f = f$ for each $j$. So, $f^2 = nf$. Now, since $n \cdot 1$ is invertible in $S$, \[ \left(\frac{f}{n\cdot 1}\right)^2 = \frac{f}{n \cdot 1}\] and the proof is complete.\end{proof}

Interesting examples of semirings such that $n \cdot 1$ is invertible for each positive integer $n$ are Gel'fand semirings. Note that a semiring $S$ is a Gel'fand semiring if $1+s \in U(S)$ for all $s\in S$ \cite[p. 56]{Golan1999(b)}. Now, it is evident that we have the following corollary:

\begin{corollary}
	Let $\{m_i\}_{i=1}^n$ be a subgroup of the monoid $M$ and $S$ a Gel'fand semiring. Then, $n^{-1} \sum_{i=1}^{n} X^{m_i}$ is a multiplicatively idempotent element of $S[M]$.  
\end{corollary}

\begin{proposition}
Let $S$ be a semiring. Then, the following statements hold:

\begin{enumerate}
	
	\item If $B$ is an $S$-semialgebra, $A$ an amount function from $B$ into $\Id(S)$, and $f$ is a $k$-potent element of $B$, then $A_f$ is a multiplicatively $k$-potent element of the semiring $\Id(S)$.
	
	\item If $M$ is a monoid and $f$ a $k$-potent element of $S[M]$, then $c(f)$ is a multiplicatively $k$-potent element of the semiring $\Id(S)$.
\end{enumerate}

\end{proposition}

\begin{proof}
(1): Let $f$ be $k$-potent. So, $f^k = f$. Since $A$ is submultiplicative, we have \[A_f = A_{f^k} \subseteq (A_f)^k \subseteq A_f.\]
(2): The proof of this point that $c$ is submultiplicative is similar to Proposition 1.1 in \cite{Rush1978} and so, omitted. Similar to what we did for the proof of the first statement, it can be proved that if $f$ is $k$-potent, then so is $c(f)$ and the proof is complete.
\end{proof}

Let us recall that a semiring $S$ is additively cancellative if $a+b = a+c$ implies $b=c$ for all $a,b,c \in S$ \cite[p. 49]{Golan1999(b)}.

\begin{theorem}
	
	\label{idempotentacsemiring}
	
	Let $S$ be an additively cancellative semiring and $X_i$ be an indeterminate over $S$ for each positive integer $i \leq k$. Then, multiplicatively idempotent elements of $S[[X_1,X_2,\ldots,X_k]]$ and $S$ are the same.
\end{theorem}
	
\begin{proof} Let $f=\sum_{n=0}^{+\infty} f_n X^n$ be multiplicatively idempotent. Clearly $f_0$ is multiplicatively idempotent. Now, let $m$ be the smallest positive integer such that $f_m$ is nonzero. Clearly, $f^2 = f$ implies that $f_0 f_m +f_0 f_m = f_m$. Since $f_0$ is multiplicatively idempotent, by multiplying the last equality by $f_0$, we obtain that $f_0 f_m +f_0 f_m = f_0 f_m$. Since $S$ is additively cancellative, we obtain that $f_0 f_m =0$. This implies that \[f_m = f_0 f_m + f_0 f_m = 0 + 0 = 0\] which is a contradiction. Consequently, multiplicatively idempotent elements of $S[[X]]$ and $S$ are the same. Finally, an obvious induction completes the proof.\end{proof}

\begin{examples}\label{idempotentacsemiringex} Theorem \ref{idempotentacsemiring} is a semiring version of Proposition 2 in \cite{KanwarLeroyMatczuk2013}. However, note that if a semiring $S$ is not additively cancellative, then Theorem \ref{idempotentacsemiring} may not hold as the following examples show:
	\begin{enumerate}
		\item  Let $a$ be a nonzero additively and multiplicatively idempotent element of a semiring $S$. Then, $f=\sum_{n=0}^{+\infty} aX^n$ is an idempotent element of $S[[X]]$.
		
		\item Let $\{f_i\}_{i=0}^{+\infty}$ be a sequence in an additively idempotent semiring $S$ such that $f_i f_j = f_{i+j}$ for all non-negative integers $i$ and $j$. Then, $f=\sum_{n=0}^{+\infty} f_n X^n$ is idempotent. For example, let $S$ be the max-plus algebra \[(\mathbb N_0 \cup \{-\infty\}, \max,+)\] and $f=\sum_{n=0}^{+\infty} n X^n$ be an element of $S[[X]]$. Now, let $e_k = k$ and imagine $c_k$ is the coefficient of $X^k$ in $f^2$. It is obvious that \[c_k = \max\{(e_0 + e_k), (e_1 + e_{k-1}), \dots, (e_{k-1} + e_1), (e_k + e_0)\} = k.\] So, it is clear that $f^2 = f$.
	\end{enumerate}
\end{examples}

\begin{theorem}
	
	\label{unitsmonoidsemiring}
	
	Let $S$ be an entire semiring and $M$ a totally ordered monoid. Then, the following statements hold:
	
	\begin{enumerate}
		\item $U(S[M]) = \{sX^m: s\in U(S),~m\in V(M)\}$.
		\item If $M$ is positive, then $U(S[M]) = U(S)$.
		\item If $M$ is positive and $S$ is a semifield, then $U(S[M]) = S\setminus\{0\}$.
	\end{enumerate}
	
	\begin{proof}
		
		(1): Let $f$ be a unit of $S[M]$ and $g$ be in $S[M]$ such that $fg=1$. By Theorem \ref{ordsemiring}, $f$ and $g$ need to be monomials. Now, let $f=sX^m$ and $g=s'X^{m'}$. It is then clear that $fg=1$ implies that $ss'=1$ and $m+m'=0$. So, $s\in U(S)$ and $m\in V(M)$.
		
		(2): It is obvious that if $s\in U(S)$, then $s\in U(S[M])$. Now, let $f\in U(S[M])$. By definition, there is a $g \in S[M]$ such that $fg =1$. This implies that $\deg(fg) = 0$. Since $S$ is entire, $\deg(f)+\deg(g) = 0$ and since $M$ is positive $\deg(f) = \deg(g) = 0$. This causes that $f\in S$ and $g\in S$. Hence, $f\in U(S)$. 
		
		(3): This is an obvious result of (2) and the proof is complete.
	\end{proof}
\end{theorem}

\begin{corollary}
	Let $S$ be an entire semiring and $X$ an indeterminate over $S$. Then, $U(S[X]) = U(S)$. In particular, if $k$ is a semifield and $X$ is an indeterminate over $k$, then $U(k[X]) = k\setminus\{0\}$. 
\end{corollary}

\begin{theorem}\label{unitspowerseriesoverzerosumfreethm}
	Let $S$ be a zerosumfree semiring. Then, \[U(S[[X_1,X_2,\ldots,X_n]]) = U(S),\] where $X_1,X_2,\ldots,X_n$ are indeterminates over $S$.
	
	\begin{proof} Let $f=\sum_{n=0}^{\infty} f_n X^n$ be a unit in $S[[X]]$. So, there is a $g=\sum_{n=0}^{\infty} g_n X^n$ in $S[[X]]$ such that $fg=1$. This implies that $f_0 g_0 =1$ and so $f_0$ (also $g_0$) is a unit in $S$. Now, since $S$ is zerosumfree, $f_i g_j = 0$ for all $i$ and $j$ except for the case $i=j=0$. In particular, $f_i g_0 =0$ for all $i \geq 1$. This implies that $f_i=0$ for all $i \geq 1$. So, $U(S[[X]]) = U(S)$. Now, by induction on the number of indeterminates over $S$, the desired statement is obtained.\end{proof}
\end{theorem}

An element $a$ of a semiring $S$ is small in $S$ if for each non-unit $b$ in $S$, $a+b$ is also non-unit \cite[p. 77]{Golan1999(b)}. In the following, we find all small elements of $S[M]$:

\begin{theorem}\label{smallelementsmonoidsemiringsthm}
	Let $S$ be a semiring and $M$ a positive monoid. Let $f_i \in S$ and $m_i \in M$ such that $m_i < m_{i+1}$ and $m_0 = 0$. Then, $f=\sum_{i=0}^{n} f_i X^{m_i}$ is small if and only if $f_0$ is small.
\end{theorem}	
	\begin{proof} Let $f$ be small in $S[M]$ and $b$ a non-unit in $S$. Therefore, by Theorem \ref{unitsmonoidsemiring}, $b$ is a non-unit in $S[M]$, and so, $f+b$ is also a non-unit. In view of Theorem \ref{unitsmonoidsemiring}, this implies $f_0 + b$ to be a non-unit in $S$. Thus, we have already proved that $f_0$ is small.
	
	Conversely, let $f_0$ be small. Take  $f=\sum_{i=0}^{n} f_i X^{m_i}$, where $f_i \in S$, $m_i \in M$, $m_i < m_{i+1}$, and $m_0 = 0$. If $g \in S[M]$ is a non-unit, then $g_0$ is also, by Theorem \ref{unitsmonoidsemiring}, non-unit. Since $f_0$ is small, $f_0 + g_0$ is also non-unit. Using Theorem \ref{unitsmonoidsemiring}, this implies $f+g$ to be non-unit and the proof is complete.\end{proof}

\section{Irreducible elements}\label{sec:irreducibleelements}

Let us recall that a nonzero element $s$ in a semiring $S$ is reducible if there are two nonzero non-unit elements $x$ and $y$ in $S$ such that $s=xy$. An element $s$ is irreducible if it is not reducible \cite{Nasehpour2019}. In view of Theorem \ref{unitsmonoidsemiring}, the proof of the following is straightforward:

\begin{proposition}
	
\label{irreduciblemonoidsemifield}
	
	Let $k$ be a semifield, $M$ a positive monoid, and $f\in k[M]$ with the positive degree. Then, $f$ is irreducible in $k[M]$ if and only if it cannot be factored into two elements of $k[M]$ of lower degree.
	
\end{proposition}

\begin{proposition}
	
	\label{irreduciblepolynomialsemifield}
	
	Let $\alpha \neq 0$ and $\beta$ be elements of an information algebra $S$ and $X$ an indeterminate over $S$. Then, $f=\alpha X^2+\beta$ cannot be factored into $g = aX+b$ and $h = cX+d$ in $S[X]$, where $a$ and $c$ are nonzero in $S$.
\end{proposition}	
	\begin{proof} Let there be two polynomials $g=aX+b$ and $h=cX+d$ in $S[X]$, where $a$ and $c$ are nonzero elements of $S$ with \[\alpha X^2+\beta = gh = acX^2 + (ad+bc)X + bd.\] Since $\beta \neq 0$, we have $bd \neq 0$. On the other hand, $ad+bc=0$. Since $S$ is zerosumfree, $ad=bc=0$. However, $a$ and $c$ are nonzero and $S$ is entire. So, $b=d=0$, a contradiction with $bd\neq 0$ and the proof is complete.
	\end{proof}

\begin{remark}
Note that Proposition \ref{irreduciblepolynomialsemifield} does not show that $f=\alpha X^2+\beta$ is necessarily irreducible. For example, $f=2X^2+4$ is reducible in $\mathbb N_0[X]$, since $2X^2+4 = 2(X^2+2)$ and both the elements 2 and $X^2+2$ are nonzero and non-unit in $\mathbb N_0[X]$. However, we have the following result for zerosumfree semifields:
\end{remark}

\begin{proposition}
\label{irreduciblepolynomialsemifield2}	Let $\alpha \neq 0$ and $\beta$ be elements of a zerosumfree semifield $k$ and $X$ an indeterminate over $k$. Then, $f=\alpha X^2+\beta$ is irreducible if and only if $\beta \neq 0$.
\end{proposition}

\begin{proof}
If $\beta = 0$, then $f=\alpha X^2 = (\alpha X)(X)$ is reducible. If $\beta \neq 0$, then the only possible way that causes $f=\alpha X^2+\beta$ to be reducible in $k[X]$ is that there are two polynomials $g=aX+b$ and $h=cX+d$ in $k[X]$, where $a$ and $c$ are nonzero elements of $k$ and $f=gh$ which is, by Proposition \ref{irreduciblepolynomialsemifield}, impossible because a zerosumfree semifield is an information algebra. This completes the proof.
\end{proof}

Let $(M,+,0,<)$ be a positive monoid. We say a nonzero element $m$ of $M$ is decomposable if there are two nonzero elements $x$ and $y$ in $M$ such that $m=x+y$. A nonzero element $m$ in $M$ is indecomposable if it is not decomposable. Now, we have the following general example of irreducible elements in the monoid semifield $k[M]$:

\begin{proposition}
	
	\label{indecomposablemonoidsemifield}
	Let $m$ be an indecomposable element of a positive monoid $M$. If $k$ is a semifield and $a \neq 0$ and $b$ are elements of $k$, then $f=aX^m + b$ is irreducible in $k[M]$.
\end{proposition}
	\begin{proof}
On the contrary, assume that $f$ is reducible and there are two nonzero non-unit elements $g$ and $h$ in $k[M]$ such that $f=gh$. This implies that \[m = \deg(f) = \deg(g) + \deg(h), \text{~with~}\deg(f) \neq 0 \text{~and~} \deg(g) \neq 0,\] contradicting the hypothesis that $m$ is indecomposable. Hence, $f$ is irreducible in $k[M]$ and the proof is complete. \end{proof}

\begin{example}
	Let $M = \mathbb N_0[\sqrt{2} ] = \{a+b \sqrt{2}: a,b \in \mathbb N_0 \}$, where by $\mathbb N_0$, we mean the set of all nonnegative integers. It is clear that $M=\mathbb N_0 [\sqrt{2}]$ is a positive monoid. Now, let $(\mathbb R^{\geq 0},+,\cdot)$ be the semifield of nonnegative real numbers. It is evident that $S=\mathbb R^{\geq 0}$ is zerosumfree. Now by Proposition \ref{indecomposablemonoidsemifield}, the elements $X+1$ and $X^{\sqrt{2}} + 1$ are irreducible in the monoid semiring $S[M] = \mathbb R^{\geq 0}[\mathbb N_0[\sqrt{2}]]$.
\end{example}

Let us recall that a polynomial $f \in S[X]$ is primitive if $c(f) = S$. The ring version (see p. 228 in \cite{Matsumura1989} and Proposition 5.17 in \cite{Aluffi2009}) of the following is due to Ferdinand Gotthold Max Eisenstein (1823--1852), a student of Ernst Eduard Kummer (1810--1893):

\begin{theorem}[Eisenstein's Irreducibility Criterion] \label{Eisensteincriterion}
	Let $S$ be a semiring and $\mathfrak{p}$ a subtractive prime ideal of $S$. Let \[f= a_n X^n + a_{n-1} X^{n-1} + \cdots +a_0\] be a polynomial in $S[X]$ such that the following conditions hold:
	
	\begin{enumerate}
		\item \label{E1} $a_n \notin \mathfrak{p}$,
		\item \label{E2} $a_i \in \mathfrak{p}$, for all $i <n$,
		\item \label{E3} $a_0 \notin \mathfrak{p}^2$. 
	\end{enumerate}
	
	Then, $f$ cannot be factored into the multiplication of non-constant polynomials in $S[X]$. Moreover, if $f$ is primitive, then $f$ is irreducible in $S[X]$.
	
	\begin{proof}
		
		Assume that $f$ can be factored into the multiplication of the following non-constant polynomials: \[ g = b_r x^r + \cdots + b_0 \text{~and~} c_s x^s + \cdots + c_0,\] where $0 <r,s<n$, and $b_r$ and $c_s$ are nonzero. Then, either $b_0$ or $c_0$ is not in $\mathfrak{p}$, since if $b_0$ and $c_0$ are both in $\mathfrak{p}$, then $a_0 = b_0 c_0$ is an element of $\mathfrak{p}^2$, contradicting the condition (\ref{E3}). Without loss of generality, we may suppose that $b_0 \notin \mathfrak{p}$. Since $a_0 \in \mathfrak{p}$ and $\mathfrak{p}$ is prime, we have $c_0 \in \mathfrak{p}$. On the other hand, since $a_n = b_r c_s$, by condition (\ref{E1}), none of the elements $b_r$ and $c_s$ are in $\mathfrak{p}$. Define $m$ to be the following number: \[m = \min \{k : c_k \notin \mathfrak{p}\}.\] Observe that
		
		\begin{displaymath}
			a_m = b_0 c_m + b_1 c_{m-1} + \cdots + \left\{ \begin{array}{ll}
				b_m c_0 & \textrm{if $r \geq m,$}\\
				b_r c_{m-r} & \textrm{if $r<m.$}\\
			\end{array} \right.
		\end{displaymath}
		
		The facts that neither $b_0$ nor $c_m$ are elements of $\mathfrak{p}$, while $c_{m-1}, \dots, c_0$ are all in $\mathfrak{p}$, and $\mathfrak{p}$ is subtractive imply that $a_m \notin \mathfrak{p}$. Therefore, $m=n$. Consequently, $s = n$, contradicting our assumption that $s<n$.
		
		Finally, let $f$ be primitive, i.e. $c(f) = S$. In order to show that $f$ is irreducible, we need to show that $f$ cannot be factored into a multiplication of a constant and a polynomial. On the contrary, assume that $f = d k$ where $d$ is a non-unit element of $S$ and $k$ is a polynomial in $S[X]$. Observe that $S = c(f) = d c(k) \subseteq \langle d \rangle$, contracting that $d$ is a non-unit. Hence, $f$ is irreducible and the proof of Eisenstein's criterion is complete.
	\end{proof}
\end{theorem}

Let us recall that a semiring $S$ is a semidomain if $ab=ac$ with $a\neq 0$ will cause $b=c$, for all $a,b,c\in S$. The definition of factorial semidomains are the same as their counterparts in ring theory, i.e., factorial domains (UFDs) \cite{Nasehpour2019}. The following is a generalization of Eisenstein's criterion for unique factorization domains (see for example, p. 147 in \cite{Spindler1994}):

\begin{corollary}\label{Eisensteincriterioncor}
Let $S$ be a weak Gaussian factorial semidomain and $p$ a prime element of $S$. Let \[f= a_n X^n + a_{n-1} X^{n-1} + \cdots +a_0\] be a polynomial in $S[X]$ such that the following condition hold:
	
	\begin{enumerate}
		\item $p \nmid a_n$,
		\item $p \mid a_i$, for all $i <n$,
		\item $p^2 \nmid a_0$. 
	\end{enumerate}
	
	Then, $f$ cannot be factored into non-constant polynomials in $S[X]$.
\end{corollary}

\begin{remark} In the following, we give some examples satisfying the suitable conditions in Corollary \ref{Eisensteincriterioncor}. Let $D$ be a Dedekind domain. By Proposition 3.10 in \cite{NasehpourV}, the semiring $\Id(D)$ with the standard addition and multiplication of ideals \cite[Proposition 6.29]{Golan1999(b)} is a factorial semidomain such that each ideal of the semiring $\Id(D)$ is subtractive. 
\end{remark}

\begin{example}
	Let $S=\Id(\mathbb Z)$ and $X$ be an indeterminate over $S$. Since $\mathbb Z$ is a Dedekind domain, $\Id(\mathbb Z)$ is a subtractive factorial semiring. Consider the following polynomial: \[f=(1) X^7 + (10)X^3 + (5)\] Observe that the principal ideal $(5)$ is a prime element of  $\Id(\mathbb Z)$ satisfying the conditions (1), (2), and (3) in Corollary \ref{Eisensteincriterioncor}. Also, note that the leading coefficient of $f$ is $(1)$, and so, $f$ is a primitive polynomial. From all we explained, we deduce that $f$ is irreducible in $S[X]$.
\end{example}

\section*{Acknowledgments}

 The paper's presentation benefited from the referee's insightful comments, for which the author is grateful.

\bibliographystyle{plain}

\end{document}